\documentclass[12pt]{article}

\usepackage{amsmath,amssymb}

\newcommand{\R}{\mathbb{R}}
\newcommand{\Z}{\mathbb{Z}}
\newcommand{\y}{{\boldsymbol{y}}}
\newtheorem{theorem}{Theorem}
\newtheorem{lemma}[theorem]{Lemma}
\newtheorem{proposition}[theorem]{Proposition}
\newtheorem{conjecture}[theorem]{Conjecture}
\newtheorem{corollary}[theorem]{Corollary}

\DeclareMathOperator{\card}{card}
\begin{document}
\title{A Sharp Inequality for Conditional Distribution of the First
Exit Time of Brownian Motion}
\date{September 7, 2004}
\author{Majid Hosseini\\
Department of Mathematics,\\
State University of New York at New Paltz,\\
75 S. Manheim Blvd. Suite 9,\\
New Paltz, NY 12561\\
hosseinm@newpaltz.edu} \maketitle \pagebreak  
\begin{abstract}
Let $U$ be a domain, convex in $x$ and symmetric about the
$y$-axis, which is contained in a centered and oriented rectangle
$R$. \linebreak If $\tau_A$ is the first exit time of Brownian
motion from $A$ and $A^+=A\cap \{(x,y):x>0\}$, it is proved that
$P^z(\tau_{U^+}>s\mid \tau_{R^+}>t)\leq P^z(\tau_{U}>s\mid
\tau_{R}>t)$ for every $s,t>0$ and every $z\in U^+$. 
\end{abstract}


\section{Introduction}

Let $A$ be a subset of  $\R^2$. The set $A$ is convex in $x$ if
its intersection with every line parallel to the $x$-axis is a
single interval or empty. We put $A^+ = A \cap \{(x,y)\mid x>0\}$.
Also, let $B_t=(B_{1,t},B_{2,t})$, $t\geq 0$ be standard two
dimensional Brownian motion  and $\tau_A=\inf\{t>0:B_t\not\in
A\}$. We will prove the following.
\begin{theorem}\label{th:compare probabilities} Let $U$ be an open, bounded, and connected set in $\R^2$
 which is symmetric about the $y$-axis
and convex in $x$. Also, let $R$ be an open rectangle containing
$U$, that is symmetric with respect to the $y$-axis, and has sides
parallel to the axes. If $z$ is a point in $U^+$, then for every
$s,t>0$,
\begin{equation}\label{eq:compare probabilities}
P^z\left(\tau_{U^+}>s\mid \tau_{R^+}>t\right)\leq
P^z\left(\tau_{U}>s\mid \tau_{R}>t\right).
\end{equation}
\end{theorem}

Recently, inequalities of this type, wherein the values of various
quantities related to
 $U$, $U^+$, $R$, and $R^+$ are compared, have been studied
extensively. Davis \cite{davis} proved the first inequality of
this kind for the heat kernel of Laplacian. Ba\~{n}uelos and
M\'{e}ndez-Hern\'{a}ndez \cite{banmen} extended Davis's result to
the heat kernel of Schr\"{o}dinger operators and integrals of
these kernels. You \cite{you} proved an inequality of this type
for the trace of Schr\"{o}dinger operators. Davis and Hosseini
\cite{DH} proved the extension to the heat content. The
inequalities studied in \cite{banmen,davis,DH,you} are ``ratio
inequalities'', in the sense that the left side of the inequality
is the ratio of some functional of $U^+$ and the same functional
for $R^+$, and the right side is the the corresponding ratio for
$U$ and $R$. Inequality~(\ref{eq:compare probabilities}) is not
strictly a ratio inequality, but rather a ``ratio-like''
inequality.

The proof of Theorem~\ref{th:compare probabilities} is based on
the idea of conditioning on zeros. Since the zeros of Brownian
motion are uncountable, it is not possible to use this approach
for Brownian motion directly. Therefore, we first prove a discrete
analog of (\ref{eq:compare probabilities}), and then use scaling.
The discrete analog is stated and proved in the following section.
Its proof is an application and modification of the techniques
introduced in \cite{DH}. We will repeat some of the material in
\cite{DH} so that we can refer to them and modify them for our
purpose. We will point them out as we go through the proof.
 The method of deriving (\ref{eq:compare probabilities})
from its discrete counterpart is a standard application of the
invariance principle. We will omit this derivation for the sake of
brevity. See \cite{DH} for a detailed description of an almost
identical derivation.

 The following example shows  that if the convexity condition in
Theorem~\ref{th:compare probabilities} is removed, we can find a
domain $U$ for which Theorem~\ref{th:compare probabilities} fails.
Let $0<d<1/2$ and take
$U=(-1,1)\times(-1,1)\setminus\{(0,y):|y|\geq d/2\}$. Let
$R=(-1,1)\times(-1,1)$ and $z=(d,1/2)$. Also, put $s=t=1$. We will
show at the end of Section~2 that, for this example, the right
side of (\ref{eq:compare probabilities}) converges to zero as
$d\rightarrow 0$ while the left side is equal to 1 for all values
of $d$.

Note that since
\begin{equation*}
E^z\left(\tau_{U}\mid \tau_{R}>t\right)= \int_0^\infty
P^z\left(\tau_{U}>s\mid \tau_{R}>t\right)\,ds,
\end{equation*}
and a similar equation holds for $E^z\left(\tau_{U^+}\mid
\tau_{R^+}>t\right)$, the following is an immediate consequence of
Theorem~\ref{th:compare probabilities}.
\begin{corollary}\label{cor: conditional expecation} Let $U$, $R$, and $z$ be as in Theorem~\ref{th:compare probabilities}.
 Then
for every $t>0$,
\begin{equation}\label{eq:conditional exopectation}
E^z\left(\tau_{U^+}\mid \tau_{R^+}>t\right)\leq
E^z\left(\tau_{U}\mid \tau_{R}>t\right).
\end{equation}
\end{corollary}

Corollary~\ref{cor: conditional expecation} and the ratio
inequalities proved in \cite{banmen,davis,DH,you} lead to the
following conjecture.

\begin{conjecture}
Let $U$, $R$, and $z$ be as in Theorem~\ref{th:compare
probabilities}. Then
\begin{equation}
\frac{E^z\left(\tau_{U^+}\right)}{E^z\left(\tau_{R^+}\right)} \leq
\frac{E^z\left(\tau_{U}\right)}{E^z\left(\tau_{R}\right)}.
\end{equation}
\end{conjecture}

\section{Discrete-Time Inequalities}
 Let $\{X_i\}_{i\geq 0}$ and
$\{Y_i\}_{i\geq 0}$ be independent sequences of random variables
such that both sequences  $\{ X_{i+1} -X_i\}$ and $\{Y_{i+1} -
Y_i\}$ are i.i.d.\ sequences of random variables, each taking
values $0$, $1$, and $-1$ with probability $1/3$. Let $Z_i = (X_i,
Y_i)$. Thus $Z_i$ is a random walk on $\Z^2$ started at $Z_0$.
Consider  ${\Lambda} \subset \Z^2$ and let  $z = (x,y)$ be a point
in $\Lambda$. For any set $A$ let $\tau_A$ be the first exit time
of $Z_i,i\geq 0$ from $A$.

\begin{proposition}\label{pr:compare probabilities discrete}
Let ${\Lambda}$ be a bounded and connected subset of $\Z^2$ which
is symmetric about the $y$-axis and convex in $x$. Let $T$ be a
 rectangle containing ${\Lambda}$ with sides parallel to the axes. Then for all
nonnegative integers $m$ and $n$,
\begin{equation}\label{eq:compare probabilities discrete}
P^z\left(\tau_{\Lambda^+}>m\mid \tau_{T^+}>n\right)\leq
P^z\left(\tau_{\Lambda}>m\mid \tau_{T}>n\right).
\end{equation}
\end{proposition}

We will prove the equivalent statement, that  for all nonnegative
integers $m$ and $n$,
\begin{equation}\label{eq:ratio of probabilities discrete}
\frac{P^z\left(\tau_{\Lambda^+}>m,
\tau_{T^+}>n\right)}{P^z\left(\tau_{\Lambda}>m,\tau_T>n\right)}
\leq
\frac{P^z\left(\tau_{T^+}>n\right)}{P^z\left(\tau_T>n\right)}.
\end{equation}

Let $l=\max(m,n)$ and let $\y = (y_0,\ldots,y_l)$ be a sequence
such that $y_0 =y$, and $|y_i - y_{i-1}| \leq 1$ for $1 \leq i
\leq l$. In addition, assume that
    \[P^z\left(Y_i = y_i , Z_i \in {\Lambda}, 0 \leq i \leq l \right)>0.\]
 We will call such $\y$ admissible. Let $\y$ be an admissible
 sequence. For an event $A$, define
 \[P^{\y,z}\left(A\right)=P^z\left(A\mid Y_i=y_i, 0\leq i\leq
 l\right).\]
 We will prove the following.
\begin{lemma}\label{lem:ratio of probabilities discrete
conditioned} Let $\Lambda$, $T$, $m$, $n$ and $z$ be as in
Proposition \ref{pr:compare probabilities discrete}. Let $\y$ be
an admissible sequence. Then
\begin{equation}\label{eq:ratio of probabilities discrete conditioned}
\frac{P^{\y,z}\left(\tau_{\Lambda^+}>m, \tau_{T^+}>n
\right)}{P^{\y,z}\left(\tau_{\Lambda}>m,\tau_T>n \right)} \leq
\frac{P^{\y,z}\left(\tau_{T^+}>n \right)}{P^{\y,z}\left(\tau_T>n
\right)}.
\end{equation}
\end{lemma}

Note that the right side of (\ref{eq:ratio of probabilities
discrete conditioned}) is independent of $\y$ and equals the right
side of (\ref{eq:ratio of probabilities discrete}). Therefore
Lemma \ref{lem:ratio of probabilities discrete conditioned}
implies Proposition \ref{pr:compare probabilities discrete}.

 For an admissible $\y$, let ${\gamma}^{\Lambda}_{\y}(k) =
\max\{i:(i,y_k) \in {\Lambda}\}$. Since $\{X_i\}_{i\geq 0}$ and
$\{Y_i\}_{i\geq 0}$ are independent, if $z=(x,y)$, then
\[P^{\y,z}\left(\tau_{T^+}>n
\right)=P^x\left(0<X_i\leq \gamma^T_{\y}(i), 0\leq i \leq
n\right),\] and
\[P^{\y,z}\left(\tau_{T}>n
\right)=P^x\left(|X_i|\leq \gamma^T_{\y}(i), 0\leq i \leq
n\right),\] and similar equalities hold for the remaining
quantities in (\ref{eq:ratio of probabilities discrete
conditioned}).

Let $f \colon \{0,\ldots,n\} \rightarrow \Z^+$ and $h \colon
\{0,\ldots,m\} \rightarrow \Z^+$ be such that $h(i)\leq f(i)$ for
all $i\in\{0,\ldots,\min(m,n)\}$. Furthermore assume $x$ is an
integer such that $x\leq h(0)$. Also, put $F^+=\{0<X_i\leq f(i),
0\leq i \leq n\}$ and $F=\{|X_i|\leq f(i), 0\leq i \leq n\}$.
Define $H^+$ and $H$ similarly, by replacing $f$ with $h$ and $n$
with $m$. We will show that
\begin{equation}\label{eq:F and H}
P^x\left(H^+\mid F^+\right) \leq P^x\left(H\mid F\right),
\end{equation}
which, by the above discussion, implies (\ref{eq:ratio of
probabilities discrete conditioned}).

In order to prove (\ref{eq:F and H}), we need to investigate the
properties of the joint distribution of $\{X_i\}_{i=0}^m$ and
$\{|X_i|\}_{i=0}^m$, given $F^+$ and $F$ respectively. For the
case $m\leq n$, this was done in \cite{DH}. In particular, Lemmas
7 and 8 in \cite{DH} are the main tools for this task. These two
Lemmas are restated as Lemmas~\ref{lem:coupling} and
\ref{lem:coupling R} in this note. Their proofs and the discussion
immediately preceding them are repeated from \cite{DH}, both for
completeness and for the modifications we will make in these
proofs to prove the case $m>n$.

 Let $a$,
$b$, $\alpha$,  and $\beta$ be integers such that  $0 \leq a < b
\leq n$ and  that $0 <{\alpha} \leq f(a)$, and $ 0 < {\beta} \leq
f(b)$. Define the probability measures $P_{a,b}^{{\alpha}}$ and
$R_{a,b}^{{\alpha},{\beta}}$ on the set of all finite sequences
$(x_a,x_{a+1},\ldots,x_b)$ of integers by
    \begin{align*}
    & P_{a,b}^{{\alpha}}(x_a,x_{a+1},\ldots,x_b) \notag\\
     & \overset{ \text{def} }{ = }P( X_k = x_k, a \leq k \leq b |
    X_a = {\alpha},
    0< X_i \leq f(i), a \leq i \leq b),\\
\intertext{and}
& R_{a,b}^{{\alpha},{\beta}}(x_a,x_{a+1},\ldots,x_b) \notag\\
     & \overset{ \text{def} }{ = }P( X_k = x_k, a \leq k \leq b \mid
    X_a = {\alpha}, X_b = {\beta},
    0< X_i \leq f(i), a \leq i \leq b).
    \end{align*}

Let $\pi_j$ be the coordinate maps:
$\pi_j(x_a,x_{a+1},\ldots,x_b)= x_j$. Under $P_{a,b}^{{\alpha}}$
the finite sequence of random variables
$\pi_a,\pi_{a+1},\ldots,\pi_b$ is a Markov chain started at
${\alpha}$ with (non-stationary) transition probabilities, which
do not depend on ${\alpha}$, given by
    \begin{equation}\label{eq:transition probability p}
    P_{a,b}^{{\alpha}}( \pi_{k+1} = v | \pi_k =u)
    =
    \frac{p_v^{k+1}}{\sum_{j=u-1}^{u+1} p_j^{k+1}} \qquad
    \mbox{for } v=u-1,u,u+1,
    \end{equation}
where
    \[ p_j^{k+1} = P(
    0 < X_i \leq f(i), k+1 \leq i \leq b
    |X_{k+1} =j).\]
Also, under $R_{a,b}^{{\alpha},{\beta}}$ the sequence
$\pi_a,\pi_{a+1},\ldots,\pi_b$ is a Markov chain started at
${\alpha}$ with  transition probabilities given by
    \begin{equation}\label{eq:transition probability R}
R_{a,b}^{{\alpha},{\beta}}( \pi_{k+1} = v | \pi_k =u)
    =
    \frac{r_v^{k+1}}{\sum_{j=u-1}^{u+1} r_j^{k+1}} \qquad
    \mbox{for } v=u-1,u,u+1,
    \end{equation}
where
    \[ r_j^{k+1} = P(X_b = {\beta},
    0 < X_i \leq f(i), k+1 \leq i \leq b
    |X_{k+1} =j).\]

\begin{lemma}\label{lem:coupling}
Let $ 0 < \alpha_0 \leq \alpha_1$. Consider Markov chains
$\psi_a(=\alpha_0),\psi_{a+1},\ldots,\psi_b$ and
$\zeta_a(=\alpha_1),\zeta_{a+1},\ldots,\zeta_b$ that have
transition probabilities given by (\ref{eq:transition probability
p}) with $\alpha$ replaced by $\alpha_0$ and $\alpha_1$
respectively.
 Then there are  Markov chains $\tilde{\psi}_a,\ldots,\tilde{\psi}_b$
and $\tilde{\zeta}_a,\ldots,\tilde{\zeta}_b$, defined on a common
probability space ${\Omega}$, and having the same transition
probabilities and initial distributions as $\psi_a,\ldots,\psi_b$
and $\zeta_a,\ldots,\zeta_b$ respectively, such that
    \begin{equation}\label{eq:coupling equation}
    \tilde{\psi}_i(\omega) \leq \tilde{\zeta}_i(\omega)
    \qquad \mbox{for every
    $\omega  \in \Omega$ and $a \leq i \leq b$}.
    \end{equation}
\end{lemma}
{\it Proof~}  We will use induction. For each $k$ satisfying $a
\leq k <  b$, all $r \in \R$, and each $u$ such that $1 \leq u\leq
f(k)$,~let~ $F_{k+1}(r,u) = P_{a,b}^{\alpha_0} ( \pi_{k+1} \leq r|
\pi_k = u)$ and $G_{k+1}(r,u) = P_{a,b}^{\alpha_1} ( \pi_{k+1}
\leq r | \pi_k = u)$ . Let $\tilde{\psi}_a = \alpha_0$ and
$\tilde{\zeta}_a = \alpha_1$ on $\Omega$. Assume that
$\tilde{\psi}_a,\ldots,\tilde{\psi}_k$ and
$\tilde{\zeta}_a,\ldots,\tilde{\zeta}_k$ have been defined on
$\Omega$. Let $T_{k+1}$ be a random variable uniformly distributed
over $\left[0,1\right]$, defined on $\Omega$, and independent of
$\tilde{\psi}_a,\ldots,\tilde{\psi}_k$ and
$\tilde{\zeta}_a,\ldots,\tilde{\zeta}_k$. For each $u$ such that
$1\leq u \leq f(k)$ and all $t \in [0,1]$ define
$F_{k+1}^{-1}(t,u) = \inf \{r : F_{k+1}(r,u) \geq t\}$ and define
$G_{k+1}^{-1}(t,u)$ similarly. On the event $\{\tilde{\psi}_k = u
\}$ let $\tilde{\psi}_{k+1} = F_{k+1}^{-1}(T_{k+1},u)$. Define
$\tilde{\zeta}_{k+1}$ in a similar manner. It is routine to check
that these are indeed Markov chains with the desired transition
probabilities and initial distributions. We use induction to prove
$\tilde{\psi}_{k}({\omega}) \leq \tilde{\zeta}_{k}({\omega})$ for
all ${\omega} \in {\Omega}$. We have $\tilde{\psi}_{a} =
{\alpha}_{0} \leq {\alpha}_{1} = \tilde{\zeta}_{a}$. Now assume
that  $\tilde{\psi}_k({\omega}) \leq \tilde{\zeta}_k ({\omega})$
for all ${\omega} \in {\Omega}$. Consider $1 \leq u \leq u'$ and
${\omega} \in  \{ \tilde{\psi}_k =u, \tilde{\zeta}_k = u'\}$. For
all $r\in \R$ and each $1 \leq u \leq u'$ we have
\begin{equation}\label{eq:stochastic dominance}F_{k+1}(r,u) \geq
G_{k+1}(r,u').
\end{equation}For if
$u'\geq u+1$, (\ref{eq:stochastic dominance}) follows directly
from the form of the transition probabilities for $\zeta$ and
$\psi$. When $u=u'$ we have $F_{k+1}(r,u)=G_{k+1}(r,u')$.  Hence
by (\ref{eq:stochastic dominance}), we have $F_{k+1}^{-1}(t,u)
\leq G_{k+1}^{-1}(t,u')$ for each $t \in [0,1]$, and so
$\tilde{\psi}_{k+1}({\omega}) \leq \tilde{\zeta}_{k+1}({\omega})$.
This completes the proof of  Lemma \ref{lem:coupling}.

\begin{lemma}\label{lem:coupling R}
Let $ 0 < \alpha_0 \leq \alpha_1$ and $ 0 < \beta_0 \leq \beta_1
$. Consider Markov chains $\psi_a (=\alpha_0),\ldots,\psi_b$ and
$\zeta_a(=\alpha_1),\ldots,\zeta_b$ that have transition
probabilities given by (\ref{eq:transition probability R}) with
$(\alpha,\beta)$ replaced by  $(\alpha_0,\beta_0)$ and
$(\alpha_1,\beta_1)$ respectively. Then there are  Markov chains
$\tilde{\psi}_a,\ldots,\tilde{\psi}_b$ and
$\tilde{\zeta}_a,\ldots,\tilde{\zeta}_b$, defined on a common
probability space ${\Omega}$, and having the same transition
probabilities and initial distributions as $\psi_a,\ldots,\psi_b$
and $\zeta_a,\ldots,\zeta_b$ respectively, such that
    \begin{equation}\label{eq:coupling equation R}
    \tilde{\psi}_i(\omega) \leq \tilde{\zeta}_i(\omega)
    \qquad \mbox{for every
    $\omega  \in \Omega$ and $a \leq i \leq b$}.
    \end{equation}
\end{lemma}
{\it Proof~}   For each $k$ satisfying $a \leq k  <  b$, all $r
\in \R$, and each $u$ such that $1 \leq u \leq f(k)$,~let~
$F_{k+1}(r,u) = R_{a,b}^{\alpha_0,\beta_0} ( \pi_{k+1} \leq r |
\pi_k = u)$ and $G_{k+1}(r,u) = R_{a,b}^{\alpha_1,\beta_1} (
\pi_{k+1} \leq r | \pi_k = u)$ . The rest of the proof follows the
proof of Lemma~\ref{lem:coupling} closely. The only difference is
in the proof of the statement that  for all $r\in \R$ and each $1
\leq u \leq u'$ we have
\begin{equation}\label{eq:stochastic dominance
R} F_{k+1}(r,u) \geq G_{k+1}(r,u').
\end{equation}
To show this first assume that $\beta_0=\beta_1$ and $u'=u+1$.
Then (\ref{eq:stochastic dominance R}) follows directly from the
form of the transition probabilities for $\zeta$ and $\psi$. When
$u=u'$ and $\beta_0=\beta_1$, we have
$F_{k+1}(r,u)=G_{k+1}(r,u')$. These facts also imply that  for
fixed $r$, the function $G_{k+1}(r,u)$ is decreasing in $u$. Put
together, the special case $\beta_0=\beta_1$ follows. Therefore,
Lemma \ref{lem:coupling R} holds in this case.  The special case
$\alpha_0=\alpha_1$ and $u'=u$ follows from the fact that running
our conditioned walks backwards in time still gives a conditioned
walk (we are just counting paths) and the special case of Lemma
\ref{lem:coupling R} for $\beta_0=\beta_1$. The general case
follows from these two special cases by first considering the
pairs $(\alpha_0,\beta_0)$ and $(\alpha_1,\beta_0)$ followed by
the pairs $(\alpha_1,\beta_0)$ and $(\alpha_1,\beta_1)$. This
completes the proof of
 (\ref{eq:stochastic dominance R}), and therefore, of Lemma
\ref{lem:coupling R}.

We can now prove (\ref{eq:F and H}).  We consider two cases: (1)
$m\leq n$ and (2) $m>n$. Let $e_0$ be the left side of (\ref{eq:F
and H}). In both cases, we will find a partition $\Psi$ of $F$
such that
    \begin{equation}\label{ineq:partition}
     e_0 \leq P^x\left(
    H\mid A\right), \quad A \in \Psi.
    \end{equation}

For the case $m\leq n$ the proof is essentially the same as the
proof of the case $m=n$ which was done in \cite{DH}. Since we will
use the ideas in this proof for proving the case $m>n$, we bring
this proof here with the minor adjustments that are necessary.

We start with the case $m\leq n$. Let $N=\card\{i: 0<i\leq n,
X_i=0\}$ and and let $M_1, M_2,\ldots,M_N$ be the indices $i\leq
n$ such that $X_i = 0$. Define
\[Q_{2,c,d} =
\{ N=2, M_1 = c, M_2 = d\}\cap F.\]

 Let $k$ be a positive integer.
For a sequence $i_1,\ldots,i_k$ of positive integers, define
$Q_{k,i_1,\ldots,i_k}$ in the same manner as $Q_{2,c,d}$. Also,
put $Q=\{N=0\}\cap F$.

Clearly $ e_0 = P^x(H|Q)$. To make the argument easier to follow,
assume that $k=2$, $M_1=c$, and $M_2=d$ and consider $Q_{2,c,d}$.
Given $Q_{2,c,d}$, the random variables $|X_0|,\ldots,|X_n|$ form
a nonhomogeneous Markov chain $\psi_0,\ldots,\psi_n$. Note that
the Markov chain $\psi_{0},\ldots,\psi_{c-1}$ has the same
transition probabilities as the sequence
$\pi_{0},\ldots,\pi_{c-1}$ under $R^{x,1}_{0,c-1}$. Also the
Markov chain $\psi_{c+1},\ldots,\psi_{d-1}$ has the same
transition probabilities  as the sequence
$\pi_{c+1},\ldots,\pi_{d-1}$ under $R^{1,1}_{c+1,d-1}$ and the
Markov chain $\psi_{d+1},\ldots,\psi_{n}$ has the same transition
probabilities as the sequence $\pi_{d+1},\ldots,\pi_{n}$ under
$P^{1}_{d+1,n}$.

Now consider $s_1$,$s_2$,$s_3 \geq 1$ and $t_1$,$t_2$,$t_3 \geq 1$
for which it is possible to condition on  $F^+$,  $X_{c-1}=s_1$,
$X_{c}=s_2$, $X_{c+1}=s_3$, $X_{d-1}=t_1$, $X_{d}=t_2$, and
$X_{d+1}=t_3$. Fix the $s_i$ and the $t_i$ ($1\leq i\leq 3$.) Then
given the above condition the  sequence $X_0,\ldots,X_n$ is a
non-homogeneous Markov chain $\zeta_0,\ldots,\zeta_n$ such that
$\zeta_{0},\ldots,\zeta_{c-1}$ will have the same transition
probabilities as $\pi_{0},\ldots,\pi_{c-1}$ under
$R^{x,s_1}_{0,c-1}$. Also $\zeta_{c+1},\ldots,\zeta_{d-1}$ has the
same transition probabilities as $\pi_{c+1},\ldots,\pi_{d-1}$
under $R^{s_3,t_1}_{c+1,d-1}$ and $\zeta_{d+1},\ldots,\zeta_{n}$
has the same transition probabilities as
$\pi_{d+1},\ldots,\pi_{n}$ under $P^{t_3}_{d+1,n}$. For the Markov
chains $\{\psi_i\}_{i=0}^{c-1}$ and $\{\zeta_i\}_{i=0}^{c-1}$,
consider
    $\tilde{\psi}_{0},\ldots,\tilde{\psi}_{c-1}$ and
    $\tilde{\zeta}_0,\ldots,\tilde{\zeta}_{c-1}$
as constructed in Lemma~\ref{lem:coupling R}. Do the same for the
time frame $[c+1,d-1]$, making sure that $\{T_i\}_{i=c+1}^{d-1}$
considered in Lemma \ref{lem:coupling R} are independent of
$\{T_i\}_{i=1}^c$, $\{\tilde{\psi}_i\}_{i=1}^c$, and
$\{\tilde{\zeta}_i\}_{i=1}^c$. For the Markov chains
$\{\psi_i\}_{i=d+1}^{n}$ and $\{\zeta_i\}_{i=d+1}^{n}$, consider
    $\tilde{\psi}_{d+1},\ldots,\tilde{\psi}_{n}$ and
    $\tilde{\zeta}_{d+1},\ldots,\tilde{\zeta}_{n}$
as constructed in Lemma~\ref{lem:coupling}, again, making sure
that $T_{d+1},\ldots,T_n$ are independent of all previous $T_i$,
$\tilde{\psi}_i$, and $\tilde{\zeta}_i$. We also have
    $s_2= \zeta_{c} > \psi_c = 0$ and $t_2=\zeta_{d} > \psi_{d}
=0$. Define $\tilde{\zeta}_{c} = s_2$, $\tilde{\psi}_{c}=0$,
$\tilde{\zeta}_{d}=t_2$, and $\tilde{\psi}_{d}=0$ on $\Omega$.

Put ${\tilde{H}}^+=\{\tilde{\zeta}_i\leq h(i),0\leq i \leq m\}$
    and $\tilde{H}=\{\tilde{\psi}_i\leq h(i), 0\leq i \leq m\}$. Hence for
    every $\omega\in \Omega$ we have $I_{\tilde{H}^+}(\omega) \leq
I_{\tilde{H}}(\omega)$.

The Markov chain ${\tilde{\psi}}$ has the same
    distribution as the Markov chain
    $|X_0|,\ldots,|X_n|$ given $Q_{2,c,d}$. The Markov chain
    $\tilde{\zeta}$ has the same
    distribution as $X_0,\ldots,X_n$ given $F^+$ and
    $X_{c-1}=s_1$,$X_c=s_2$,$X_{c+1}=s_3$,$X_{d-1}=t_1$,$X_d=t_2$,
    and $X_{d+1}=t_3$.

 Therefore for all possible values of  $s_1$, $s_2$, $s_3$, $t_1$, $t_2$, and $t_3
\geq 1$,

    \begin{align*}
     P^x\Bigl( & H^+  \mid F^+ \mbox{ and }  X_{c-1}=s_1, X_c=s_2, X_{c+1}= s_3,
    X_{d-1}=t_1, X_d =t_2, X_{d+1} =t_3\Bigr)\\
        &\leq
    P^x\Bigl(H \mid Q_{2,{c},{d}}\Bigr).
    \end{align*}
Since this is true for all possible values of
$s_1$,$s_2$,$s_3$,$t_1$,$t_2$,$t_3 \geq 1$, we have that
\begin{equation*}  P^x\left(
H^+ \mid F^+ \right) \leq
\\ P^x\left( H \mid
Q_{2,c,d}\right).
\end{equation*}
 This proves (\ref{ineq:partition}) for $Q_{2,c,d}$. The
argument for all other $Q_{k,{i_1},\ldots,{i_k}}$ is the same as
above. Therefore the proof of (\ref{ineq:partition}) for the case
$m\leq n$ is complete.

To prove (\ref{ineq:partition}) for the case $m>n$, let $N$,
$M_1,\ldots,M_N$, $Q$ and $Q_{k,i_1,\ldots,i_k}$ be as in the
proof of the case $m\leq n$. Again, to make the proof easier to
follow, we focus on $Q_{2,c,d}$. Given $Q_{2,c,d}$, the random
variables $|X_0|,\ldots,|X_m|$ form a nonhomogeneous Markov chain
$\psi_0,\ldots,\psi_m$. The part $\psi_0,\ldots,\psi_n$ has the
same transition probabilities as discussed in the proof of the
case $m\leq n$. The part $\psi_n,\ldots,\psi_m$ has transition
probabilities
\begin{align}\label{eq:transition probability psi}
    P( \psi_{k+1} = v | \psi_k =u)
     &= \frac{1}{3}
     \qquad
    \mbox{for }  u >0 \mbox{ and }v=u-1,u,u+1
    ,\\
    \intertext{and }
P( \psi_{k+1} = 1 | \psi_k =0)
     &= \frac{2}{3}=1-P( \psi_{k+1} = 0 | \psi_k =0).
    \end{align}

Now consider $ s_1$,$s_2$,$s_3\geq 1 $ and $ t_1$,$t_2$,$t_3 \geq
1$ for which it is possible to condition on $F^+$, $X_{c-1}=s_1$,
$X_{c}=s_2$, $X_{c+1}=s_3$, $X_{d-1}=t_1$, $X_{d}=t_2$, and
$X_{d+1}=t_3$. Fix the $s_i$ and the $t_i$. Then given the above
condition the  sequence $X_0,\ldots,X_m$ is a non-homogeneous
Markov chain $\zeta_0,\ldots,\zeta_m$. The part
$\zeta_0,\ldots,\zeta_n$ has the same transition probabilities as
discussed in the proof of the case $m\leq n$. The part
$\zeta_n,\ldots,\zeta_m$ has transition probabilities
\begin{equation}\label{eq:transition probabilitiy zeta}
P(\zeta_{k+1} = v \mid \zeta_k=u) =\frac{1}{3} \qquad \mbox{for
}v=u-1,u,u+1.
\end{equation}

Construct $\tilde{\psi}_{0},\ldots,\tilde{\psi}_{n}$ and
$\tilde{\zeta}_{0},\ldots,\tilde{\zeta}_{n}$ as in the proof of
the case $m\leq n$. For $n\leq k\leq m-1$, put
$F_{k+1}(r,u)=P(\psi_{k+1}\leq r \mid \psi_k=u)$ and
$G_{k+1}(r,u)=P(\zeta_{k+1}\leq r\mid \zeta_k=u)$. It is easy to
see that for $0<u \leq v$ with $v>0$,
\begin{equation}\label{eq:stochastic dominance after n}
F_{k+1}(r,u)\geq G_{k+1}(r,v) \qquad \mbox{for } n\leq k\leq m-1.
\end{equation}

Next, we will construct
$\tilde{\psi}_{n+1},\ldots,\tilde{\psi}_{m}$ and
$\tilde{\zeta}_{n+1},\ldots,\tilde{\zeta}_{m}$ on $\Omega$ by the
method used in Lemma \ref{lem:coupling}. Assume that
$\tilde{\psi}_n,\ldots,\tilde{\psi}_k$ and
$\tilde{\zeta}_n,\ldots,\tilde{\zeta}_k$ have been defined
on~$\Omega$. Let $T_{k+1}$ be a random variable uniformly
distributed over $\left[0,1\right]$, defined on~$\Omega$,
independent of $T_0,\ldots,T_k$, as well as independent of
$\tilde{\psi}_0,\ldots,\tilde{\psi}_k$ and
$\tilde{\zeta}_0,\ldots,\tilde{\zeta}_k$. For each $u$ such that
$0\leq u $ and all $t \in [0,1]$ define $F_{k+1}^{-1}(t,u) = \inf
\{r : F_{k+1}(r,u) \geq t\}$ and for $u \in \Z$ define
$G_{k+1}^{-1}(t,u)$ similarly. On the event $\{\tilde{\psi}_k = u
\}$ let $\tilde{\psi}_{k+1} = F_{k+1}^{-1}(T_{k+1},u)$. Define
$\tilde{\zeta}_{k+1}$ in a similar manner. Again, it is easy to
check that $\tilde{\psi}_{n},\ldots,\tilde{\psi}_{m}$ and
$\tilde{\zeta}_{n},\ldots,\tilde{\zeta}_{m}$ are Markov chains
with the same transition probabilities as $\psi_n,\ldots,\psi_m$
and $\zeta_n,\ldots,\zeta_m$ respectively.

We will show that
\begin{align}\{0<\tilde{\zeta_k}\leq h(k),0\leq
k\leq m \}& \subseteq \{|\tilde{\psi_k}| \leq h(k),0\leq k\leq
m\}\\
& = \{\tilde{\psi_k}\leq h(k), 0\leq k \leq m\}.\notag
\end{align}

Note that, by construction, for $0\leq k\leq n$ and for all
$\omega \in \Omega$, we have $\tilde{\psi_k}(\omega) \leq
\tilde{\zeta_k}(\omega)$ and $0<\tilde{\zeta_k}(\omega)$. We will
show, by induction, that for $n+1\leq k\leq m$, if
$0<\tilde{\zeta_k}(\omega)$, then $\tilde{\psi_k}(\omega) \leq
\tilde{\zeta_k}(\omega)$. First note that
$\tilde{\psi_n}(\omega)\leq \tilde{\zeta_n}(\omega)$ and
$0<\tilde{\zeta_n}(\omega)$ for all $\omega \in \Omega$.

Now assume that $\omega \in \{\tilde{\psi}_k\leq
\tilde{\zeta}_k,\tilde{\zeta}_k >0\}$. We will show that
$\tilde{\psi}_{k+1}(\omega)\leq \tilde{\zeta}_{k+1}(\omega)$.
Consider integers $u$ and $v$ with $0\leq u\leq v$ and $v>0$ and
assume that $\omega\in\{\tilde{\psi}_k=u,\tilde{\zeta}_k=v\}$. By
(\ref{eq:stochastic dominance after n}), we have
$F_{k+1}^{-1}(t,u)\leq G_{k+1}^{-1}(t,v)$ for each $t\in
\left[0,1\right]$. Therefore, if $\omega \in \{0<\tilde{\zeta}_k
\leq h(k), 0\leq k\leq m\}$, then $\omega \in \{\tilde{\psi}_k\leq
h(k),0\leq k\leq m\}$. Hence for all $\omega \in \Omega$,
\[I_{\{0<\tilde{\zeta}_k\leq h(k),0\leq k\leq m\}}(\omega) \leq
I_{\{\tilde{\psi}_k\leq h(k), 0\leq k \leq m\}}(\omega).\]

The Markov chain $\tilde{\psi}_0,\ldots,\tilde{\psi}_m$ has the
same
    distribution as the Markov chain
    $|X_0|,\ldots,|X_m|$ given $Q_{2,c,d}$. The Markov chain
    $\tilde{\zeta}_0,\ldots,\tilde{\zeta}_m$ has the same
    distribution as $X_0,\ldots,X_m$ given $F^+$ and
    $X_{c-1}=s_1$,$X_c=s_2$,$X_{c+1}=s_3$,$X_{d-1}=t_1$,$X_d=t_2$,
    and $X_{d+1}=t_3$.

 Therefore for all possible values of  $s_1, s_2, s_3, t_1, t_2,t_3
\geq 1$,

    \begin{align*}
    P^x\Bigl(
     H^+\mid & F^+ \mbox{ and } X_{c-1}=s_1, X_c=s_2, &X_{c+1}= s_3,
    X_{d-1}=t_1, X_d =t_2, &X_{d+1} =t_3\Bigr)\\
        & \leq
    P^x\Bigl(H \mid Q_{2,{c},{d}}\Bigr).
    \end{align*}
Since this is true for all possible values of
$s_1,s_2,s_3,t_1,t_2,t_3 \geq 1$, we have that
\begin{equation*}  P^x\left(H^+ \mid F^+ \right) \leq
 P^x\left( H \mid Q_{2,c,d}\right).
\end{equation*}
 This proves  for $Q_{2,c,d}$. The
argument for $Q$ and all other $Q_{k,{i_1},\ldots,{i_k}}$ is the
same as above. Therefore the proof of (\ref{ineq:partition}) and
of Lemma \ref{lem:ratio of probabilities discrete conditioned} is
complete. As we pointed out earlier, Lemma \ref{lem:ratio of
probabilities discrete conditioned} implies Proposition
\ref{pr:compare probabilities discrete}. Finally, the functional
central limit theorem implies that Theorem \ref{th:compare
probabilities} follows from Proposition \ref{pr:compare
probabilities discrete}. This completes the proof of Theorem
\ref{th:compare probabilities}.

Now we will show that for the example in the introduction
inequality (\ref{eq:compare probabilities}) fails. Recall that
$B_t=(B_{1,t},B_{2,t})$, $t\geq 0$, is the standard two
dimensional Brownian motion.  Note that $R^+=U^+$ and therefore
the left side of (\ref{eq:compare probabilities}) equals 1. On the
other hand for any $\theta <1$ we have,
\begin{align}\label{eq:example}
P^z\left(\tau_U>1,\tau_R>1\right)\leq P^z\left(\tau_U>1\right)
\leq &P^z\left(B_{1,t}>0\mbox{ for all } t\in
(0,\theta)\right)\\\notag
 & + P^z\left(
B_{2,t} <d/2 \mbox{ for some } t\in (0,\theta)\right).
\end{align}

Given $\epsilon >0$, by the continuity of Brownian motion paths,
there exists $\theta_0 \in (0,1)$ such that
\begin{align}\label{eq:example continuity 1}
&P^z\left(B_{2,t} <d/2\mbox{ for some }
t\in(0,\theta_0)\right)\leq\notag\\ &P^z\left(B_{2,t}<1/4\mbox{
for some }t\in (0,\theta_0)\right) <\frac{\epsilon}{2}.
\end{align}

Also, by the Law of iterated logarithm, for $d$ small enough, we
have
\begin{equation}\label{eq:example continuity 2}
P^z\left(B_{1,t}>0 \mbox{ for all } t\in
(0,\theta_0)\right)<\frac{\epsilon}{2}.
\end{equation}

Put together, (\ref{eq:example}), (\ref{eq:example continuity 1}),
and (\ref{eq:example continuity 2}) imply that
\begin{equation}\label{eq:limit of tau U}
\lim_{d\rightarrow 0} P^z(\tau_U>1, \tau_R>1)=0.
\end{equation}

 On
the other hand,
\begin{equation}\label{eq:limit of tau R}
\lim_{d\rightarrow
0}P^z(\tau_R>1)=P^{(0,1/2)}(\tau_R>1)>0.
\end{equation}

It follows immediately from (\ref{eq:limit of tau U}) and
(\ref{eq:limit of tau R}) that
\[\lim_{d\rightarrow
0}P^z(\tau_U>1\mid\tau_R>1)=0.\]

\section{Higher Dimensions}
The analog of Theorem \ref{th:compare probabilities} holds for an
arbitrary dimension $k$. We state it as Theorem~\ref{th:Theorem 2}
and we will show how the proof of Theorem \ref{th:compare
probabilities} can be modified to prove it.

Represent a point in  $\R^k$ by $z=(z_1,\cdots,z_k)$. Also for
$A\subseteq \R^k$ put $A^+=A\cap \{z\in \R^k \mid z_1>0\}$. Call
$A$ convex in $z_1$ if the intersection of $A$ with every line
parallel to $z_1$-axis is a connected interval or empty.

\begin{theorem}\label{th:Theorem 2} Let $k$ be a positive integer
 and let $U$ be a bounded, connected, and open subset of $\R^k$
 which is symmetric about $\{z_1=0\}$
and convex in $z_1$. Also, let
$R=(-L_1,L_1)\times\ldots\times(-L_k,L_k)$ be a $k$-dimensional
rectangle, containing $U$. Then for every $z\in U^+$ and every
$s,t>0$,

\begin{equation}\label{eq:Theorem 2}
P^z\left(\tau_{U^+}>s\mid \tau_{R^+}>t\right)\leq
P^z\left(\tau_{U}>s\mid \tau_{R}>t\right).
\end{equation}
\end{theorem}
The proof of Theorem \ref{th:Theorem 2} is very similar to the
proof of Theorem~\ref{th:compare probabilities}. Let $\{X_i^1\}_{i
\geq 0}$ ,\ldots,$\{X_i^k\}_{i \geq 0} $ be $k$ independent
one-dimensional random walks, each constructed as $\{X_i\}_{i\geq
0}$ in Section~2. Then the analog of Proposition~\ref{pr:compare
probabilities discrete} holds by conditioning on
$\{(X_i^2,\cdots,X_i^k)\}_{i=0}^l$.

{\bf Acknowledgement} I am most grateful to Professor Burgess
Davis for his valuable comments on the earlier drafts of this
paper, and in particular for pointing out the example in the
introduction.

\end{document}